\documentclass[
]{article}
\usepackage{amsmath,amssymb}
\usepackage{booktabs}
\usepackage{iftex}
\usepackage{booktabs,array}
\usepackage{graphicx}
\usepackage[authoryear,round]{natbib}
\usepackage{arxiv}
\usepackage{orcidlink}
\usepackage{amsmath}
\usepackage[T1]{fontenc}
\usepackage{url}
\hypersetup{
    hidelinks
}

\newcommand{\runninghead}{A Preprint }
\renewcommand{\runninghead}{THK-AI Research Report 1/2026 }
\title{Physics-based Approximation and Prediction of Speedlines in
Compressor Performance Maps}
\def\asep{\And}
\def\arowsep{\AND}
\author{\textbf{Abdul-Malik Akiev}\\\\TH Köln,
Germany\\\\\href{mailto:abdulmalikakiev@yahoo.com}{abdulmalikakiev@yahoo.com}\asep\textbf{Danyal
Ergür}\\\\TH Köln,
Germany\\\\\href{mailto:danyal.erguer@outlook.com}{danyal.erguer@outlook.com}\arowsep\textbf{Alexander
Schirger}\\\\TH Köln,
Germany\\\\\href{mailto:alexander.schirger02@gmail.com}{alexander.schirger02@gmail.com}\asep\textbf{Matthias
Müller}\\\\Everllence SE, Turbolader
Entwicklung\\\\\href{mailto:matthias.mueller.c@everllence.com}{matthias.mueller.c@everllence.com}\arowsep\textbf{Alexander
Hinterleitner}~\orcidlink{0009-0002-7615-6952}\\\\THK-AI Research Cluster, TH
Köln\\\\\href{mailto:alexander.hinterleitner@th-koeln.de}{alexander.hinterleitner@th-koeln.de}\asep\textbf{Thomas
Bartz-Beielstein}~\orcidlink{0000-0002-5938-5158}\\\\THK-AI Research Cluster, TH
Köln\\\\\href{mailto:thomas.bartz-beielstein@th-koeln.de}{thomas.bartz-beielstein@th-koeln.de}}
\date{}
\begin{document}
\twocolumn[{%
\maketitle
\begin{abstract}
Speedlines in compressor performance maps (CPMs) are critical for
understanding and predicting compressor behavior under various operating
conditions. We investigate a physics-based method for reconstructing
compressor performance maps from sparse measurements by fitting each
speedline with a superellipse and encoding it as a compact,
interpretable vector (surge, choke, curvature, and shape parameters).
Building on the formulation of Llamas et al., we develop a robust
two-stage fitting pipeline that couples global search with local
refinement. The approach is validated on industrial data-sets for
different turbocharger types. We discuss prediction quality for inter-
and extrapolation, metric sensitivities and outline opportunities for
physics-informed constraints, alternative function families, and hybrid
physics--ML mappings to improve boundary behavior and, ultimately,
enable full CPM reconstruction from limited data.
\end{abstract}
\vspace{1em}
}]

\section{Introduction}\label{sec-introduction}

Compressor performance maps (CPMs) play a crucial role in the design,
evaluation, and control of turbocharged engines and gas turbines. These
maps characterize the relationship between the corrected mass flow, the
pressure ratio, and the rotational speed of a compressor. However,
generating accurate CPMs requires extensive measurements on test
benches, which are time-consuming, expensive, and often limited to
discrete operating points and speedlines. Therefore, there is a strong
demand for reliable algorithms that can reconstruct or predict missing
speedlines or even entire CPMs from sparse measurements \citep{shen2019}. Recent advances in machine learning have shown promising results
in approximating CPMs using data-driven methods such as recurrent neural
networks (RNNs) \citep{schulz2025}, which model the structure of
speedlines as sequences and learn their behavior from a large number of
maps. While these approaches achieve high accuracy in interpolation and
extrapolation tasks, they require significant amounts of training data
and often lack physical interpretability.
In this report, we present a slightly modified physics-based approach to
CPM reconstruction and prediction based on a geometric fitting
procedure. Specifically, we approximate each speedline using a
\emph{superellipse}, as proposed by \citet{Llamas2019}, and represent each fitted speedline via a compact
\(\beta\)-vector containing physically meaningful parameters. Based on
this representation, we interpolate and extrapolate \(\beta\)-vectors to
predict unknown speedlines at other rotational speeds. Our approach
offers three key advantages: a low-dimensional, interpretable
parameterization of speedlines, excellent interpolation capability
particularly in mid-speed regions, and minimal data requirements,
relying on only a few known speedlines.

The full Python implementation used for fitting, prediction, and evaluation
is publicly available in our GitHub repository (\href{https://github.com/sequential-parameter-optimization/spotspeedlines/tree/teamproject}{sequential-parameter-optimization/spotspeedlines}).
The repository contains the code base used for the experiments reported
in this paper and is intended to support transparency and reproducibility.

This paper is structured as follows: Section~\ref{sec-background}
describes the background and motivation for this work, including the
industrial context and related research. Section~\ref{sec-method}
details our method, including the superellipse fitting procedure, the
prediction strategy using \(\beta\)-vectors, and the evaluation
framework. Section~\ref{sec-experiments} presents the experiments to assess interpolation and extrapolation performance. Section~\ref{sec-results} discusses the related results in detail. Finally,
Section~\ref{sec-discussion} discusses the implications of our findings,
limitations of the current approach, and potential directions for future
research.

\section{Background and Motivation}\label{sec-background}

This project was conducted as part of the so-called ``Teamprojekt''
module within the Bachelor's program in Electrical
Engineering---Automation Technology at the THK-AI Research Cluster at TH
Köln (\url{https://thk-ai.de}). The module encourages students to independently design and
implement a technical-scientific project in collaboration with an
external industrial partner. The project is embedded in a joint research
initiative with Everllence (formerly known as MAN Energy Solutions), a
company specializing in high-performance turbomachinery. As part of this
collaboration, Everllence has provided access to an internal dataset
consisting of approximately 15 CPMs, measured from two different types
of turbochargers. These real-world datasets serve as the basis for the
prediction and model validations conducted in this work.

Within this research initiative, a second team is developing a machine
learning-based approach using recurrent neural networks (RNNs) to model
CPMs. Their work focuses on exploiting sequential structures within
speedlines to interpolate and extrapolate operating data using
data-driven techniques. Their results are documented in a separate
scientific paper forming the core of the machine-learning branch of this
study \citep{schulz2025}.

In contrast to the machine-learning approach, our team investigates a
physics-based alternative approach. Our method models each speedline
using a superellipse parameterization, inspired by the work of \citet{Leufven2013} and further extended by
\citet{Llamas2017}. We encode each
speedline into a compact \(\beta\)-vector consisting of interpretable
parameters (e.g., surge, choke, curvature) and then interpolate or
extrapolate between these vectors to predict unknown speedlines.

The overarching motivation of this work is exploratory in nature. We aim
to assess how well a physically structured approximation method can
reconstruct compressor maps when provided with limited measurement data.
Although full CPM prediction is not the primary goal, the long-term
aspiration is to enable robust, low-data reconstructions of entire
performance maps using physically meaningful abstractions. By aligning
our evaluation criteria with those of the RNN-based approach, we seek to
enable a fair and quantitative comparison between physics-informed and
machine learning methods in the domain of turbomachinery modeling.

\section{Method}\label{sec-method}

\subsection{Superellipse Model}\label{superellipse-model}

Our implementation builds upon the superellipse model proposed by
\citet{Llamas2019}, extending it with robust optimization
strategies. The model approximates each speedline using a parameterized
superellipse curve, characterized by a \(\beta\)-vector containing
physically meaningful parameters such as surge point, choke point, and
shape parameters. The superellipse formulation following \citet{Llamas2019} is given by:

\begin{equation}\label{eq:superellipse}
\left(\frac{\dot{m} - \dot{m}_{zs}}{\dot{m}_{ch} - \dot{m}_{zs}}\right)^{\operatorname{CUR}} + \left(\frac{\Pi - \Pi_{zs}}{\Pi_{ch} - \Pi_{zs}}\right)^{\operatorname{CUR}} = 1,
\end{equation} 
where \(\dot{m}_{zs}\) and \(\Pi_{zs}\) represent the surge point,
\(\dot{m}_{ch}\) and \(\Pi_{ch}\) the choke point, and
\(\operatorname{CUR}\) the curvature parameter that controls the shape
of the speedline. These parameters form the \(\beta\)-vector
\(\boldsymbol{\beta} = [\dot{m}_{zs},  \Pi_{zs},  \dot{m}_{ch},  \Pi_{ch},  \operatorname{CUR}]^T\)
\citep{shen2019}.

The key innovation in our approach lies in the optimization procedure.
Instead of relying on simple initialization with heuristics and local
optimization, we implement a multi-stage fitting process. The
optimization procedure consists of three main stages. First, the global
parameter initialization employs either Particle Swarm Optimization
(PSO) with configurable swarm size (default: 100 particles, 50
iterations) or Differential Evolution (DE) with population size 15 and
up to 1,000 iterations \citep{Storn1997, Ken95}. Both methods thoroughly explore the parameter
space to avoid local minima. The second stage implements a local
optimization pipeline using L-BFGS-B as the primary solver with a
maximum of 5,000 iterations, automatically falling back to Nelder-Mead
if L-BFGS-B fails \citep{neld65a,Broy70a}. This is followed by comprehensive bounds validation
to ensure solution validity.

For performance assessment, we employ multiple error metrics. While the
above description outlines the default configuration, all parameters in
both the fitting process and later prediction can be customized through
a configuration interface. This includes swarm sizes, iteration counts,
solver choices, error metrics, etc. The described settings represent our
defaults based on the first experiments, but users can adjust these
parameters to suit specific requirements or performance needs.

\subsection{\texorpdfstring{Prediction Using \(\beta\)
Vectors}{Prediction Using \textbackslash beta Vectors}}\label{predicting-using-beta-vectors}

After fitting individual speedlines to obtain \(\beta\)-vectors, we use
polynomial interpolation to predict unknown speedlines. The process is
identical regardless of whether the target speed lies within the known
range (interpolation) or outside it (extrapolation). The prediction
process follows a systematic approach with three key components. The
hold-out strategy begins by removing the target speedline from the
dataset, followed by fitting all remaining speedlines to obtain their
\(\beta\)-vectors, which are then used to predict the held-out
speedline. For polynomial fitting, we employ 4th-degree polynomials
fitted independently to each \(\beta\) parameter, with optional speed
normalization to the {[}0,1{]} range for numerical stability. Physical
constraints are enforced throughout the process, such as maintaining
curvature \(\geq 2\) and ensuring non-negative values. The final
parameter prediction stage involves evaluating these polynomials at the
target speed to obtain the predicted \(\beta\)-vector, followed by
applying post-processing constraints to ensure physical validity of the
results.

Similar to the fitting process, all prediction parameters can be
customized through a configuration class. This includes the polynomial
degree, normalization settings, parameter constraints and the evaluation
mode (pressure-based or mass flow-based). The described settings
represent our defaults based on the first experiments, but users can
adjust these parameters according to their specific needs. Our approach can be seen as an extension of the approach by \citet{shen2019}, who use a polynomial degree of 2 and a normalization to the {[}0,1{]} range.

\subsection{Evaluation and Error
Assessment}\label{sec-error-metrics}

After predicting a \(\beta\)-vector for the target speedline, we
evaluate the prediction quality using multiple error metrics. Our
evaluation framework encompasses three comprehensive aspects. For
prediction evaluation, we implement two distinct modes: a pressure-based
approach comparing predicted pressure ratios \(\pi\) at given mass flows
\(\dot{m}\), and a mass flow-based approach comparing predicted mass
flows \(\dot{m}\) at given pressure ratios \(\pi\). Both modes utilize
actual measurement points from the held-out speedline for comparison.
The assessment employs multiple error metrics.

For model assessment we employ RMSE (absolute units), MAPE (relative
\%), the residual standard deviation (solver consistency), and the
orthogonal distance \[
\operatorname{ortho} = \sum_{i=1}^{N} \left\| \mathbf{x}^{\text{true}}_{i} - \mathbf{x}^{\text{pred}}_{i} \right\|^{2}
\]

in the \((\dot{m}, \pi)\) space. Because each metric has its own scale
and interpretation, direct numerical aggregation or comparison is
inappropriate \citep{hyndman2006another}. RMSE accentuates large
outliers, MAPE scales each error by its true value (and fails when y
\(\approx\) 0), the residual SD captures run--to--run stability, whereas
\(\operatorname{ortho}\) measures the geometric fidelity of the entire
speed-line. We therefore report all metrics separately and consider a
solver successful only when it achieves simultaneously low means and
dispersions in at least two metrics.
Our statistical-analysis
methodology ensures robustness by calculating both mean errors and
standard deviations, carefully handling numerical edge cases such as
division by zero and invalid predictions, and providing comprehensive
error statistics for thorough model comparison.

\section{Experiments}\label{sec-experiments}

\subsection{Dataset and Preprocessing}\label{dataset-and-preprocessing}

Our experiments are conducted on a simulated dataset, a public benchmark
dataset and a comprehensive dataset provided by Everllence, consisting
of 15 compressor performance maps from two different turbocharger types.
Each map contains multiple speedlines with varying numbers of
measurement points across different operating conditions. The
preprocessing involves two main steps: normalization of mass flow and
pressure ratio values, and consistent speedline identification together with
grouping to ensure data quality.

\subsection{Fitting Optimization
Study}\label{fitting-optimization-study}

The first experimental phase focused on establishing the most robust
fitting procedure for individual speedlines (no prediction). Our
benchmark study evaluated three key aspects. For model implementations, we compared the
\texttt{LlamasEllipse} implementation following \citet{Llamas2019} with the \texttt{DirectEllipse} approach based on direct least squares fitting
\citep{Fitzgibbon1999}. The \texttt{DirectEllipse} method transforms
the geometric problem into a generalized eigenvalue system: \[
\mathbf{S}\mathbf{a} = \lambda\mathbf{C}\mathbf{a},
\]
where \(\mathbf{a} = [a, b, c, d, e, f]^T\) represents the coefficients
of the general conic equation 
\[
ax^2 + bxy + cy^2 + dx + ey + f = 0,
\]
and \(\mathbf{S}\) and \(\mathbf{C}\) are scatter and constraint
matrices ensuring the elliptical shape.

We investigated three optimization strategies: a baseline approach using
direct L-BFGS-B optimization, a PSO-enhanced method , and a DE-enhanced approach.
For performance evaluation, we employed the error metrics
discussed in Section~\ref{sec-error-metrics}. Each configuration was
evaluated using simulated data to establish baseline performance and
identify the most robust parameter estimation approach.

\subsection{Prediction Method
Validation}\label{prediction-method-validation}

The second phase validated our prediction methodology using a public
benchmark dataset (\texttt{tca\_88}): Our validation methodology
comprised three key aspects. The hold-out Strategy involved systematic
testing of different speedline hold-out positions, comprehensive
analysis of interpolation versus extrapolation performance, and thorough
validation of polynomial degree selection. For parameter constraints, we
focused on implementing physical bounds for \(\beta\)-vectors,
validating surge/choke point predictions, and assessing the stability of
the curvature parameter. The numerical stability analysis encompassed
testing of speed normalization effects, evaluating polynomial fitting
robustness, and developing comprehensive prediction confidence metrics.

\subsection{Industrial Dataset
Evaluation}\label{industrial-dataset-evaluation}

The final experimental phase applied our optimized approach to the
complete Everllence dataset: The evaluation consisted of two main
components. Our cross-validation approach implemented leave-one-out
validation for each speedline, coupled with detailed statistical
analysis of prediction errors and a thorough comparison of interpolation
versus extrapolation accuracy. The performance assessment phase
evaluated the method across two different turbocharger types, each with
distinct characteristics in terms of speed ranges, number of speedlines,
and measurement point density. This comprehensive evaluation allowed us
to quantify prediction reliability across diverse operational
conditions.

\section{Results}\label{sec-results}

\subsection{Fitting Optimization
Study}\label{fitting-optimization-study-1}

Our comprehensive benchmark of different fitting strategies revealed
several key insights about the optimization of the superellipse
parameters (Eq.~\ref{eq:superellipse}). The analysis focused on comparing different
combinations of initialization strategies, solvers, and error metrics.

\begin{figure*}
\centering{
\includegraphics[width=1.0\linewidth,height=\textheight,keepaspectratio]{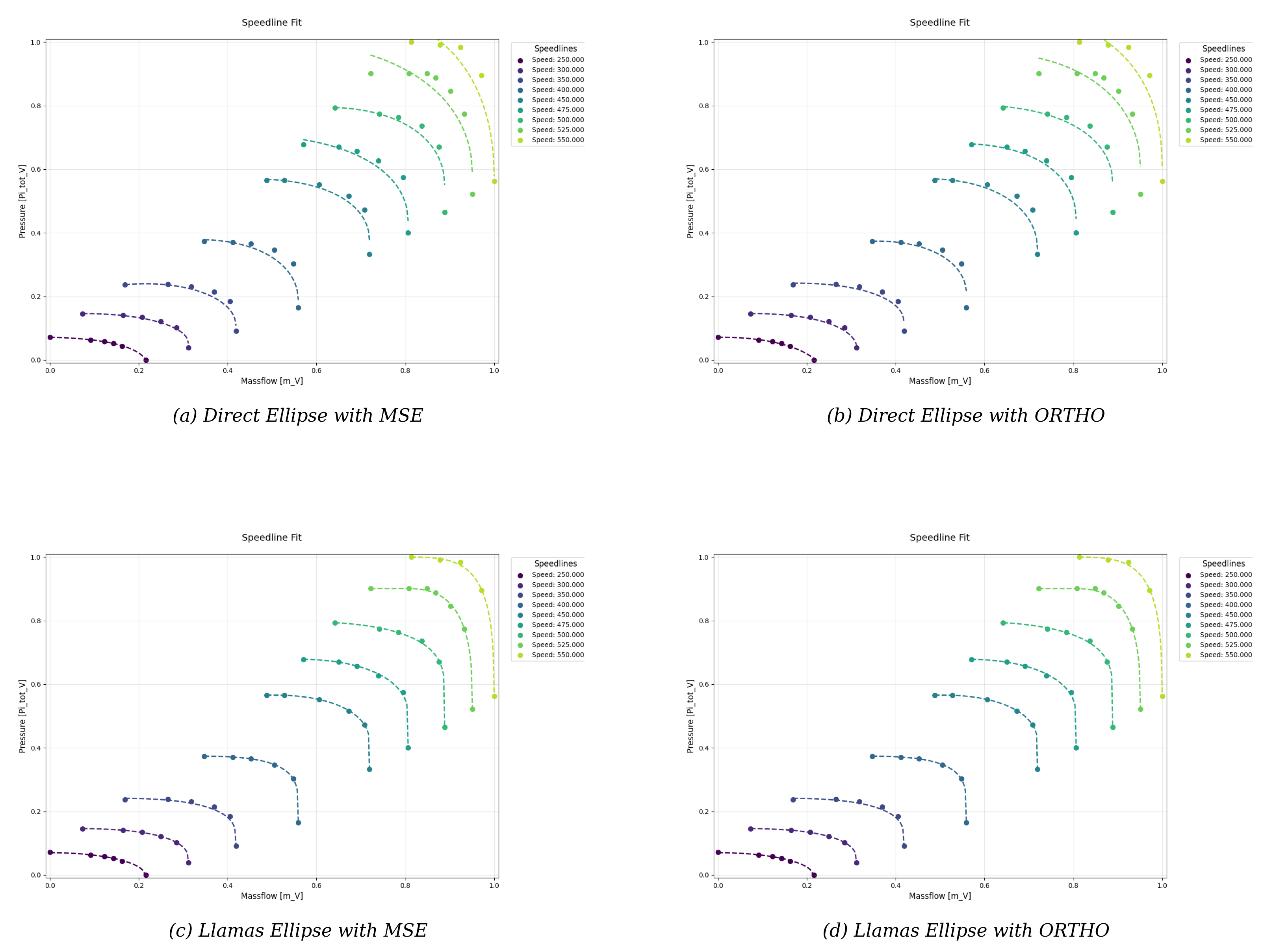}
}
\caption{\label{fig-comparison-modelperformance}Comparison of
\texttt{DirectEllipse} and \texttt{LlamasEllipse} models using the best performing
configuration on the tca88 dataset}
\end{figure*}

The \texttt{LlamasEllipse} model consistently outperformed the \texttt{DirectEllipse}
approach across multiple metrics. The model achieved lower values across
all error metrics, though this alone does not guarantee better
performance. More importantly, it demonstrated superior performance in
capturing the overall shape of the speedline, which is crucial for
accurate compressor modeling.
Figure~\ref{fig-comparison-modelperformance} compares the performance of
the \texttt{LlamasEllipse} and \texttt{DirectEllipse} models using the best performing
configuration on the \texttt{tca88} dataset, demonstrating the superior
performance of the \texttt{LlamasEllipse} model across all error metrics.

The orthogonal distance (\(\operatorname{ortho}\)) metric emerged as the
most suitable objective function. It excels in several aspects: it
provides better geometric accuracy of fitted curves, shows enhanced
robustness against outliers compared to MSE/RMSE/MAPE, and maintains
consistent performance across different speed ranges. Most importantly,
it serves as an intuitive geometric metric that is both easy to
interpret and numerically stable, avoiding NaN or infinite values in
calculations.

The combination of DE with Nelder-Mead as solver and
\(\operatorname{ortho}\) as error metric proved to be the most robust on
the \texttt{LlamasEllipse} model. Importantly the initial guess optimization
proved crucial, drastically enhancing final optimization performance and
reducing maximum error by 67\(\%\) compared to direct optimization. DE
optimization of initial parameters slightly outperformed PSO, achieving
a mean \(\operatorname{ortho}\) improvement of 84\(\%\) compared to
81\(\%\). Additionally, Nelder-Mead demonstrated superior convergence
characteristics compared to L-BFGS-B, particularly for challenging
speedline configurations.

\begin{figure*}
\centering{
\includegraphics[width=1.0\linewidth,height=\textheight,keepaspectratio]{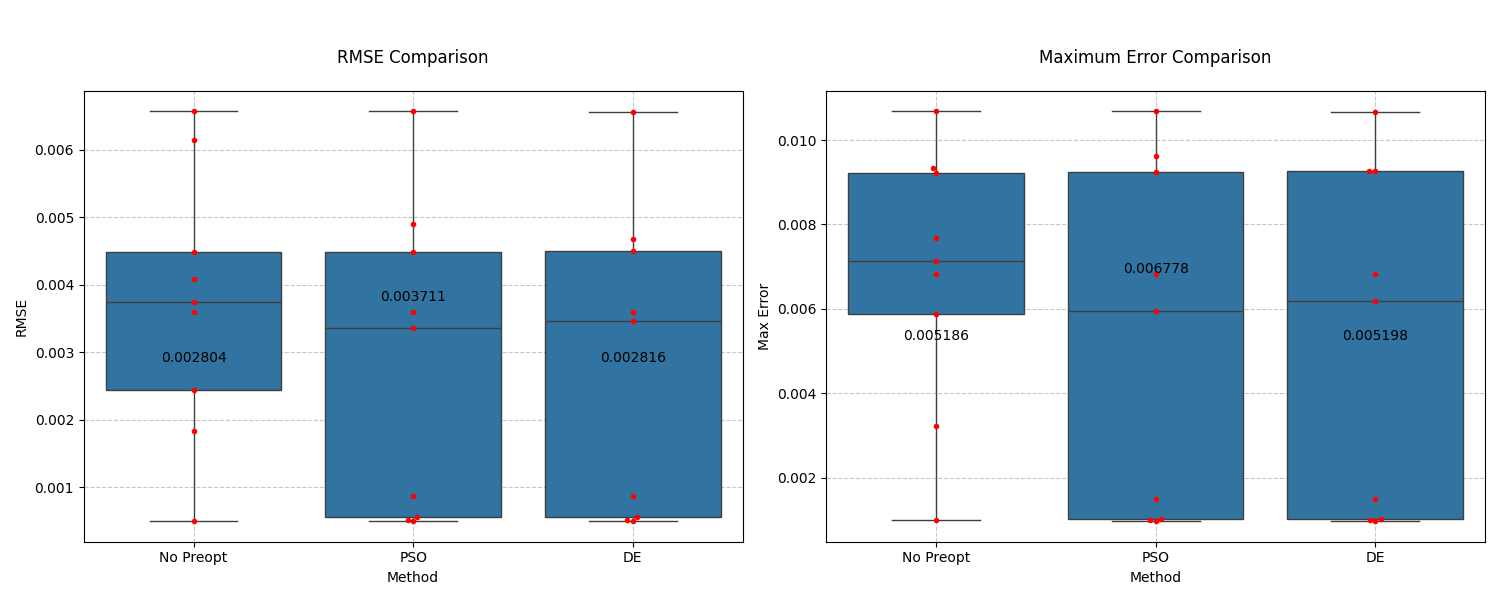}
}
\caption{\label{fig-comparison-optimization}Comparison of optimization
strategies. The boxplots show the distribution of RMSE and maximum error
across all speedlines in the tca88 dataset. Individual data points (red)
show the actual error values for each speedline.}
\end{figure*}

Based on these results, we selected the \texttt{LlamasEllipse} model with DE
initialization (population size 15, 1000 iterations), Nelder-Mead as
solver, and \(\operatorname{ortho}\) as the error metric for all
subsequent experiments. This configuration provided the best balance
between fitting accuracy and numerical stability.
Figure~\ref{fig-comparison-optimization} shows a detailed comparison of
the three optimization approaches (no pre-optimization, PSO, and DE) in
terms of RMSE and maximum error on the \texttt{tca88} dataset. The
results show that the pre-optimization provides better results compared
to the direct optimization approach.

\subsection{Interpolation Within CPMs}\label{interpolation-within-cpms}

\begin{figure*}
\centering{
\includegraphics[width=0.95\linewidth,height=\textheight,keepaspectratio]{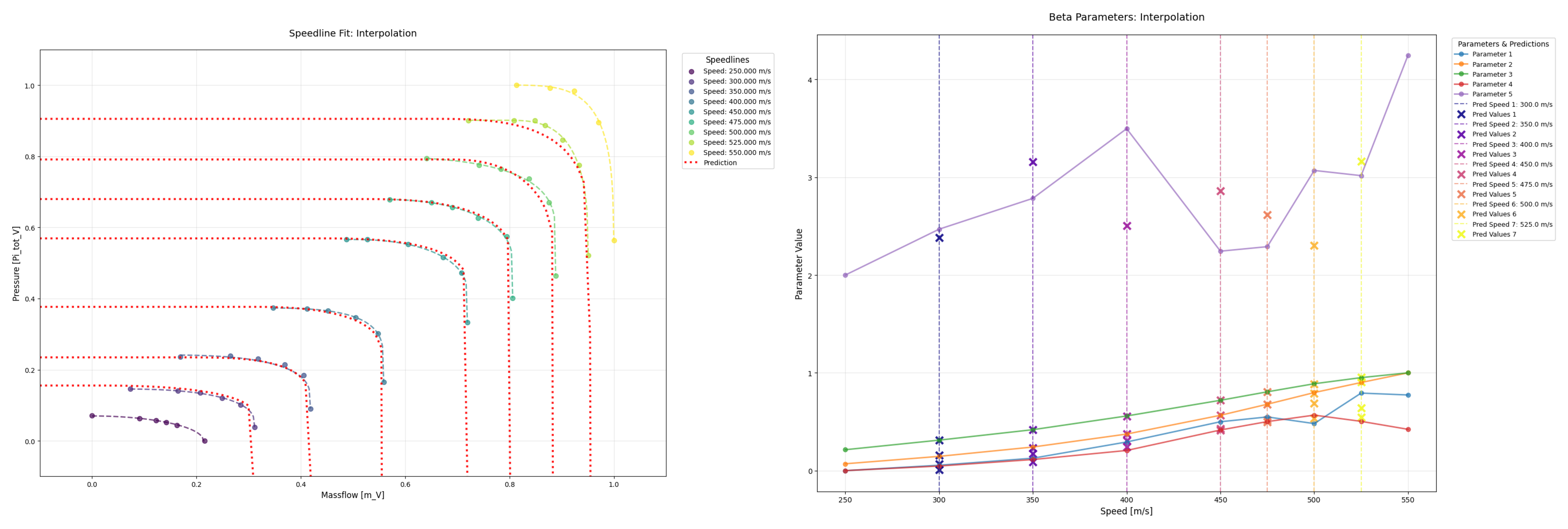}
}
\caption{\label{fig-interpolation-tca88}Interpolation results for
\texttt{tcaA88} compressor. Left: Visualization of predicted speedlines
(dashed) compared to actual data (solid) for speeds within the training
range. Right: Evolution of model coefficients (\(\beta\) vectors) across
different speeds, showing the polynomial interpolation used for
prediction.}
\end{figure*}

The interpolation results for the \texttt{tca88} dataset demonstrate the
effectiveness of our prediction approach, as shown in
Figure~\ref{fig-interpolation-tca88}. The left plot illustrates the
predicted speedlines overlaid with the actual measurements, while the
right plot shows the evolution of the \(\beta\) coefficients across
different speeds, providing insight into the interpolation behavior.

A notable observation from our experiments is that MAPE proves to be an
unreliable metric for assessing prediction quality in our context. For
instance, at 400 m/s (Index 3), we observe a MAPE of
\(161.23\% \pm 356.93\%\), which suggests poor performance, yet the
visual fit is reasonable. This discrepancy arises because MAPE is highly
sensitive to small absolute differences when the reference values are
close to zero, which is common in our normalized data. Instead, the
orthogonal distance metric \(\operatorname{ortho}\) provides more
meaningful insights into prediction quality. Most predictions show
excellent \(\operatorname{ortho}\) values below 0.01, indicating high
accuracy. The best results were achieved at 475 m/s (Index 5) with
\(\operatorname{ortho} = 0.0013 \pm 0.0027\), followed by 450 m/s (Index
4) with \(\operatorname{ortho} = 0.0031 \pm 0.0057\), and 525 m/s (Index
7) with \(\operatorname{ortho} = 0.0059 \pm 0.0070\).

An interesting case is the prediction at 350 m/s (Index 2), which shows
a relatively high \(\operatorname{ortho}\) value of
\(2.0656 \pm 4.0088\) despite visually excellent fit quality, same
applies for the speed 500 m/s (Index 6). This apparent contradiction can
be attributed to the orthogonal distance being more sensitive to small
deviations in regions where the speedline has high curvature,
particularly near the choke and surge points. The RMSE of
\(0.0165 \pm 0.0005\) for this case better reflects the actual
prediction quality we observe visually. The RMSE generally aligns well
with visual assessment, with most predictions achieving values below
0.06. Two notable exceptions occur for speed at 400 m/s with
\(\operatorname{RMSE} = 0.6456 \pm 0.9319\) and speed at 500 m/s with
\(\operatorname{RMSE} = 1.1437 \pm 2.9239\)

\begin{table}
\centering
\tiny
\setlength{\tabcolsep}{4pt} 
\caption{Metrics for speedline predictions (interpolation) on the TCA88 dataset. Each metric is shown with its standard deviation.}\label{tbl-interpolation-metrics}
\begin{tabular}{lrrrr}
\toprule
Index & Speed [m/s] & RMSE & MAPE [\%] & ORTHO \\
\midrule
1 & $300.00$ & $0.06 \pm 0.01$ & $64.51 \pm 129.75$ & $0.01 \pm 0.01$ \\
2 & $350.00$ & $0.02 \pm 0.00$ & $9.44 \pm 13.97$ & $2.07 \pm 4.01$ \\
3 & $400.00$ & $0.65 \pm 0.93$ & $161.23 \pm 356.93$ & $0.42 \pm 0.93$ \\
4 & $450.00$ & $0.05 \pm 0.01$ & $7.38 \pm 13.64$ & $0.00 \pm 0.01$ \\
5 & $475.00$ & $0.04 \pm 0.00$ & $3.98 \pm 7.74$ & $0.00 \pm 0.00$ \\
6 & $500.00$ & $1.14 \pm 2.92$ & $102.37 \pm 224.04$ & $2.40 \pm 3.41$ \\
7 & $525.00$ & $0.06 \pm 0.01$ & $5.80 \pm 9.02$ & $0.01 \pm 0.01$ \\
\bottomrule
\end{tabular}
\end{table}

These higher RMSE values occur at speeds where the speedline shape
changes more dramatically, making the interpolation more challenging.
However, as seen in the left plot of
Figure~\ref{fig-interpolation-tca88}, even these predictions capture the
general trend of the speedline behavior.
Table~\ref{tbl-interpolation-metrics} provides a comprehensive overview
of all prediction metrics for each speedline.

\subsection{Extrapolation beyond CPMs}\label{extrapolation-beyond-cpms}

\begin{figure*}
\centering{
\includegraphics[width=1.0\linewidth,height=\textheight,keepaspectratio]{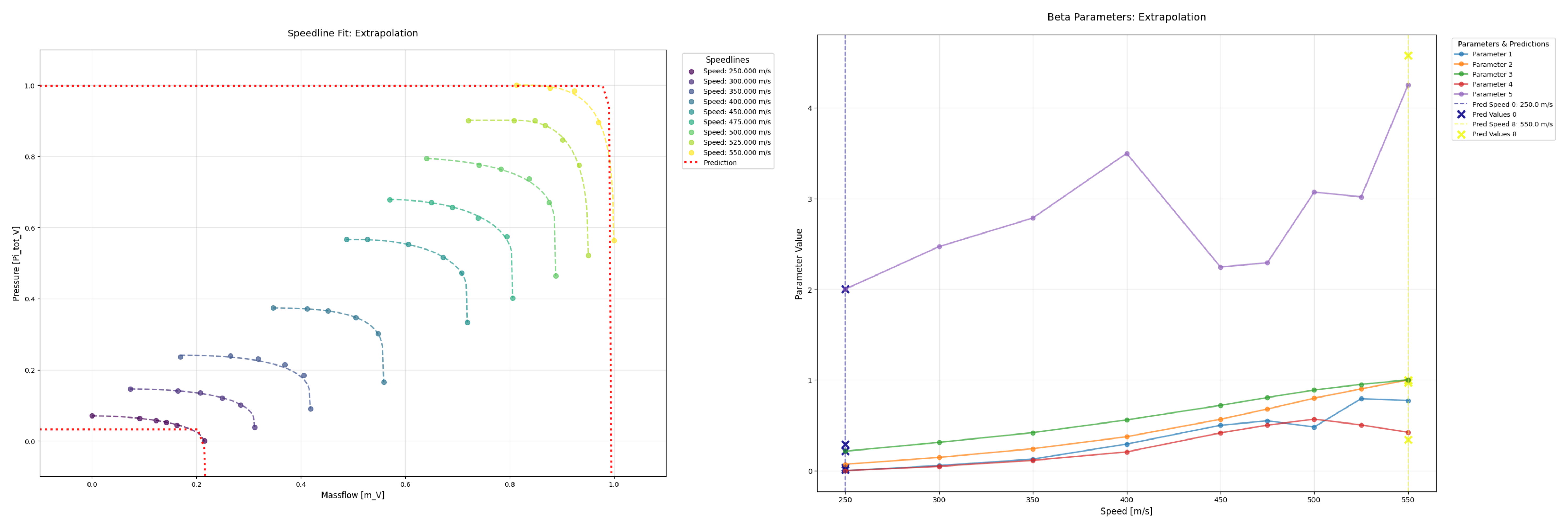}
}
\caption{\label{fig-extrapolation-tca88}Extrapolation results for TCA88
compressor. Left: Visualization of predicted speedlines (dashed)
compared to actual data (solid) for speeds outside the training range.
Right: Evolution of model coefficients (\(\beta\) vectors) across
different speeds, showing the polynomial interpolation used for
prediction.}
\end{figure*}

The extrapolation of speedlines beyond the range of available training
data presents significantly greater challenges compared to
interpolation. Despite the seemingly logical progression of the
\(\beta\) coefficients across speeds (as shown in the right plot of
Figure~\ref{fig-extrapolation-tca88}), the prediction quality
deteriorates substantially when attempting to predict speedlines at the
boundaries of the compressor map.

Table~\ref{tbl-extrapolation-metrics} presents the
complete metrics for the extrapolation predictions.
The most striking example is the prediction at 250 m/s (Index 0), where
we observe catastrophically high error metrics:
\(\operatorname{ortho} = 48194.0615 \pm 41051.4939\) and an astronomical
MAPE of over 1,000,000\%. This extreme deviation clearly indicates that
the polynomial interpolation approach, while effective for
interpolation, fails to capture the physical behavior of the compressor
at lower speeds. The prediction at 550 m/s (Index 8) performs relatively
better but still shows considerable error with
\(\operatorname{RMSE} = 0.6387 \pm 0.8107\) and
\(\operatorname{ortho} = 1.4420 \pm 2.0146\).

While mathematically sound, the polynomial interpolation of \(\beta\)
coefficients appears unable to capture the underlying physical
relationships governing compressor behavior at extreme speeds. The model
shows asymmetric reliability in extrapolation, with significantly worse
performance at lower speeds (250 m/s) compared to higher speeds (550
m/s). In addition to these findings, the smooth progression of \(\beta\)
coefficients proves to be a potentially misleading indicator of
prediction quality, highlighting the critical importance of thorough
validation for extrapolation results.

\begin{table}
\centering
\tiny
\setlength{\tabcolsep}{4pt} 
\caption{Metrics for speedline predictions (extrapolation) on the TCA88 dataset. Each metric is shown with its standard deviation.}\label{tbl-extrapolation-metrics}
\begin{tabular}{lrrrr}
\toprule
Index & Speed [m/s] & RMSE & MAPE [\%] & ORTHO \\
\midrule
0 & $250.00$ & $0.04 \pm 0.00$ & $1.02\mathrm{e}{6} \pm 2.27\mathrm{e}{6}$ & $4.82\mathrm{e}{4} \pm 4.11\mathrm{e}{4}$ \\
8 & $550.00$ & $0.64 \pm 0.81$ & $53.31 \pm 99.90$ & $1.44 \pm 2.01$ \\
\bottomrule
\end{tabular}
\end{table}
These results strongly suggest that alternative approaches should be
considered for extrapolation tasks. While the current method excels at
interpolation within the measured speed range, its limitations become
evident when pushed beyond these boundaries. This finding has important
implications for compressor map modeling and suggests that
physics-informed constraints or different mathematical models might be
necessary for reliable extrapolation.

\subsection{Performance on the Industrial
Dataset}\label{performance-on-the-industrial-dataset}

In this section we evaluate our approach on the complete industrial
dataset comprising two turbocharger types and 15 CPMs. Again, we
distinguish between interpolation and extrapolation. Overall,
interpolation delivers promising visual fidelity, whereas extrapolation
is unreliable. Error metrics are only partially informative due to
extreme outliers and metric-specific pathologies (especially MAPE), so
we complement quantitative results with visual inspection.

\subsubsection{Reliability of Metrics and Effect of
Outliers}\label{reliability-of-metrics-and-effect-of-outliers}

A few catastrophic cases inflate averages by several orders of
magnitude, especially for Type-1 maps. Means should therefore be
interpreted with caution. Robust summaries (medians, trimmed means)
would give a more balanced view. This indicates that outliers dominate
means. Regarding the sensitivity of the metrics, we observed that MAPE
becomes unstable when denominators approach zero in normalized data,
resulting in deceptively large values despite good visual alignment.
\(\operatorname{ortho}\) captures geometric fidelity well but can spike
near high-curvature regions (surge/choke) even when overall fit looks
acceptable. \(\operatorname{RMSE}\) often aligns best with visual
impression but is still affected by local shape changes.

The very high Type‑1 averages are caused by a handful of extreme
failures. Many Type‑1 predictions nonetheless look visually plausible,
as summarized in Table~\ref{tbl-type-summary}.

\begin{table}
\centering
\tiny
\setlength{\tabcolsep}{4pt} 
\caption{Summary statistics comparing interpolation and extrapolation performance across both compressor types. Note the dramatic increase in error metrics highlighting the method's limitations beyond the training range.}\label{tbl-type-summary}
\begin{tabular}{lrrrr}
\toprule
Metric & Type-0 (Int.) & Type-1 (Int.) & Type-0 (Ext.) & Type-1 (Ext.) \\
\midrule
Number of CPMs & 6 & 9 & 6 & 9 \\
Total predictions & 37 & 38 & 12 & 18 \\
Average RMSE & $0.69$ & $30.15$ & $2.12\mathrm{e}{3}$ & $41.85$ \\
Average MAPE [\%] & $94.27$ & $2.53\mathrm{e}{4}$ & $7.10\mathrm{e}{10}$ & $5.85\mathrm{e}{8}$ \\
Average ORTHO & $4.23$ & $3.16\mathrm{e}{4}$ & $5.26\mathrm{e}{8}$ & $3.63\mathrm{e}{8}$ \\
\bottomrule
\end{tabular}
\end{table}

\subsubsection{Industrial Dataset: Interpolation}\label{interpolation-across-cpms-promising-visual-quality}

Interpolation works well across most CPMs of the industrial dataset. Despite occasional metric
spikes, predicted speedlines usually track the ground truth closely,
especially in mid-speed regions. Lower speedlines perform worse than
higher speedlines. Predictions are more reliable for CPMs with more
speedlines and more measurement points.

We observed three typical situations across CPMs. Results can be visually solid, although metrics are relatively poor due to extreme outliers in specific indices.
This is precisely where pure numeric results become misleading: the plot
shows the predicted curve adhering to the measured shape for large
portions of the speedline. Sometimes situations with very
accurate interpolation results, with the predicted speedline closely
following the measured curve across the entire range, can occur. Both visual
inspection and error metrics confirm high fidelity, highlighting the
robustness of the superellipse fitting in well-constrained mid-speed
regions. Finally, interpolation can be very poor, especially for low-speed lines with fewer measurement points.

Visual inspection exemplifies the divergence between raw metrics and
perceived geometric fidelity. It demonstrates that visually satisfactory
predictions can co-exist with discouraging metrics, underlining the
importance of qualitative assessment. CPMs can be acceptable in shape despite large reported errors, e.g., if isolated poor indices but an overall
reasonable fit occurs.

\subsubsection{Industrial Dataset: Extrapolation}\label{extrapolation-beyond-cpms-unreliable}

Extrapolation (first/last speedlines) is also on the industrial dataset substantially less reliable.
Error distributions are heavy‑tailed with catastrophic outliers (e.g.,
type‑0 and type‑1 boundary indices), indicating that polynomial
\(\beta\)‑interpolation is not physics‑grounded enough at the borders.
Visuals often show meaningful trends, but quantitative errors remain
high. For industrial application, we recommend restricting to
interpolation unless additional physics-based constraints or alternative
models are introduced for boundary behavior.

\section{Discussion}\label{sec-discussion}

The results of the fitting optimization study clearly demonstrate that
the combination of the \texttt{LlamasEllipse} model, DE initialization,
Nelder--Mead local solver, and the orthogonal distance metric provides
the most reliable and consistent results across different datasets. Each
of these elements contributes a specific advantage: The \texttt{LlamasEllipse}
model inherently captures the geometric characteristics of compressor
speedlines more effectively than the \texttt{DirectEllipse} approach. While
\texttt{DirectEllipse} minimizes algebraic error, it often struggles to represent
the curved structure of real speedlines, especially near surge and choke
regions. In contrast, the \texttt{LlamasEllipse} formulation directly
parameterizes physically meaningful points, resulting in better fidelity
and interpretability.

The optimization process benefits greatly from global search strategies.
DE proved particularly effective at providing robust initial guesses,
avoiding local minima that hindered the performance of direct or
PSO-based initialization. The stochastic population-based nature of DE
offers a more thorough exploration of the parameter space, which is
crucial given the non-convexity of the fitting problem. For local
convergence, Nelder--Mead outperformed L-BFGS-B in difficult cases.
While L-BFGS-B is efficient for smooth, well-conditioned landscapes, it
can fail in the presence of irregularities or poorly scaled gradients,
both common in compressor data. Nelder--Mead, being derivative-free,
converged more robustly and provided stable final fits, especially when
combined with high-quality DE initialization. The orthogonal distance
metric emerged as the most suitable objective function. Unlike MSE,
RMSE, or MAPE, \(\operatorname{ortho}\) evaluates the geometric distance
between fitted and measured points, aligning directly with the physical
goal of reproducing the overall shape of the speedline. Its robustness
against outliers and its avoidance of division-by-zero issues make it a
stable and interpretable choice. Taken together, these choices explain
why the selected configuration (\texttt{LlamasEllipse} + DE + Nelder--Mead +
\(\operatorname{ortho}\)) yielded the most robust performance in the
fitting study. It combines physical interpretability, robust global
exploration, stable local refinement, and a geometrically meaningful
error metric.

When moving from fitting to prediction, the results show a clear
asymmetry between interpolation and extrapolation. Interpolation within
the measured speed range works reliably: predicted speedlines generally
follow the measured curves closely, with low RMSE and
\(\operatorname{ortho}\) values, particularly in mid-speed regions where
\(\beta\)-vectors evolve smoothly. This indicates that polynomial
interpolation of \(\beta\)-parameters is a valid strategy as long as the
target speed lies within the convex hull of the known data.

By contrast, extrapolation beyond the available range proved unreliable.
Although \(\beta\)-vectors show a smooth progression across speeds,
polynomial models fail to capture the true physical behavior at the
boundaries of compressor maps. This is particularly evident at lower
speeds, where catastrophic errors occur despite apparently plausible
\(\beta\)-trajectories. The mismatch highlights a key limitation:
polynomial extrapolation is mathematically consistent but not
physics-grounded. It can therefore produce misleading predictions if not
carefully validated.
While the presented approach demonstrates strong interpolation
capabilities and robust fitting performance, several important
limitations remain. The most evident limitation lies in the poor
extrapolation performance. Polynomial interpolation of
\(\beta\)-vectors, while effective within the measured speed range,
fails to generalize beyond the boundaries of the compressor map. At low
speeds in particular, predictions diverge drastically despite smooth
coefficient trajectories. This underscores the fact that polynomial
regression, being purely mathematical, does not encode the underlying
physics of compressor behavior. Comparable findings have been reported
by \citet{shen2019}, who also noted that low-degree
polynomials provide satisfactory interpolation but are unreliable for
extrapolation.

Sparse data reduces the constraints available for \(\beta\)-vector fitting and weakens
polynomial stability, as observed in maps such as CPM (0,4). This
limitation aligns with observations in previous work by \citet{Llamas2017}, where higher data density was
identified as a key enabler for stable parameterization.
Although \(\operatorname{ortho}\) provides a geometrically meaningful
measure, it can overemphasize local deviations at surge and choke
points, leading to seemingly poor numerical scores despite visually
convincing fits. Conversely, MAPE often inflates errors near small
denominators in normalized data. This suggests that relying on a single
metric is insufficient. A combination of RMSE, \(\operatorname{ortho}\),
and visual inspection should be maintained for comprehensive evaluation.

Several opportunities exist to improve prediction quality and extend
applicability. Incorporating physical constraints (e.g., monotonic
behavior near surge, bounds derived from thermodynamic limits) could
guide extrapolation and prevent catastrophic divergence. Hybrid
strategies that combine polynomial regression with physics-based
boundary conditions represent a promising direction. Instead of
polynomial regression, spline interpolation or rational functions could
provide more stable coefficient evolution across speeds, reducing
oscillations and improving extrapolation robustness. \citet{shen2019} employed quadratic fits; however, more flexible functions
may strike a better balance between bias and variance.
A natural extension is the combination of interpretable \(\beta\)-vector
parameterization with machine learning methods. Neural networks or
Gaussian processes could be trained to map rotational speed to
\(\beta\)-vectors, offering data-driven flexibility while retaining a
physics-informed representation. This would bridge the gap between
purely physics-based and purely data-driven methods, similar to the
RNN-based approach of \citet{schulz2025}.

In summary, while the proposed method achieves reliable interpolation
with minimal data requirements, it remains limited at the edges of the
compressor map and under sparse conditions. Future work should address
extrapolation through physics-informed constraints, alternative function
families, and hybrid ML integration. Such extensions could eventually
enable robust reconstruction of entire CPMs, aligning with the long-term
vision of data-efficient, interpretable compressor modeling.

\section*{Acknowledgments}\label{sec-acknowledgments}

This research and development project was funded by the German Federal
Ministry of Education and Research (Bundesministerium für Bildung und
Forschung, BMBF) within the funding measure \emph{Forschung an
Fachhochschulen -- KI-Nachwuchs@FH 2-2021} under the project
\emph{TH Köln -- Künstliche Intelligenz plus (THK-KIplus)},
funding code 13FH007KI2. The authors are responsible for the content of
this publication.

\bibliographystyle{plainnat}
\bibliography{references}

\begin{thebibliography}{11}
\providecommand{\natexlab}[1]{#1}
\providecommand{\url}[1]{\texttt{#1}}
\expandafter\ifx\csname urlstyle\endcsname\relax
  \providecommand{\doi}[1]{doi: #1}\else
  \providecommand{\doi}{doi: \begingroup \urlstyle{rm}\Url}\fi

\bibitem[Broyden(1970)]{Broy70a}
C~G Broyden.
\newblock {The Convergence of a Class of Double-Rank Minimization Algorithms}.
\newblock \emph{Journal of the Institute of Mathematics and Its Applications},
  6:\penalty0 76--90, 1970.

\bibitem[Fitzgibbon et~al.(1999)Fitzgibbon, Pilu, and Fisher]{Fitzgibbon1999}
Andrew Fitzgibbon, Maurizio Pilu, and Robert~B. Fisher.
\newblock Direct least squares fitting of ellipses.
\newblock \emph{IEEE Transactions on Pattern Analysis and Machine
  Intelligence}, 21\penalty0 (5):\penalty0 476--480, 1999.
\newblock \doi{10.1109/34.765658}.

\bibitem[Hyndman and Koehler(2006)]{hyndman2006another}
Rob~J. Hyndman and Anne~B. Koehler.
\newblock Another look at measures of forecast accuracy.
\newblock \emph{International Journal of Forecasting}, 22\penalty0
  (4):\penalty0 679--688, 2006.
\newblock \doi{10.1016/j.ijforecast.2006.03.001}.

\bibitem[Kennedy and Eberhart(1995)]{Ken95}
J~Kennedy and R~Eberhart.
\newblock {Particle Swarm Optimization}.
\newblock In \emph{Proceedings of the Fourth IEEE International Conference on
  Neural Networks}, pages 1942--1948, Piscataway NJ, 1995. IEEE.

\bibitem[Leufv{\'e}n and Eriksson(2013)]{Leufven2013}
Oskar Leufv{\'e}n and Lars Eriksson.
\newblock A surge-and-choke capable compressor flow model -- validation and
  extrapolation capability.
\newblock \emph{Control Engineering Practice}, 21\penalty0 (12):\penalty0
  1871--1883, 2013.
\newblock \doi{10.1016/j.conengprac.2013.07.005}.

\bibitem[Llamas and Eriksson(2017)]{Llamas2017}
Xavier Llamas and Lars Eriksson.
\newblock Parameterizing compact and extensible compressor models using
  orthogonal distance minimization.
\newblock \emph{Journal of Engineering for Gas Turbines and Power},
  139\penalty0 (1):\penalty0 012601--1--012601--10, 2017.
\newblock \doi{10.1115/1.4034152}.

\bibitem[Llamas and Eriksson(2019)]{Llamas2019}
Xavier Llamas and Lars Eriksson.
\newblock Control-oriented modeling of two-stroke diesel engines with exhaust
  gas recirculation for marine applications.
\newblock \emph{Proceedings of the Institution of Mechanical Engineers, Part M:
  Journal of Engineering for the Maritime Environment}, 233\penalty0
  (2):\penalty0 551--574, 2019.

\bibitem[Nelder and Mead(1965)]{neld65a}
J.~A. Nelder and R.~Mead.
\newblock {A Simplex Method for Function Minimization}.
\newblock \emph{The Computer Journal}, 7\penalty0 (4):\penalty0 308--313, 01
  1965.
\newblock ISSN 0010-4620.
\newblock \doi{10.1093/comjnl/7.4.308}.
\newblock URL \url{https://doi.org/10.1093/comjnl/7.4.308}.

\bibitem[Schulz et~al.(2025)]{schulz2025}
R.~Schulz et~al.
\newblock Automated prediction of compressor performance maps:surrogate-based
  optimization with rnns for enhanced extrapolation and interpolation.
\newblock \emph{arXiv preprint arXiv:2506.12345}, 2025.

\bibitem[Shen et~al.(2019)Shen, Zhang, Zhang, Yang, and Jia]{shen2019}
Haosheng Shen, Chuan Zhang, Jundong Zhang, Baicheng Yang, and Baozhu Jia.
\newblock Applicable and comparative research of compressor mass flow rate and
  isentropic efficiency empirical models to marine large-scale compressor.
\newblock \emph{Energies}, 12\penalty0 (24):\penalty0 4749, 2019.
\newblock \doi{10.3390/en13010047}.
\newblock URL \url{https://www.mdpi.com/1996-1073/13/1/4749}.

\bibitem[Storn and Price(1997)]{Storn1997}
Rainer Storn and Kenneth Price.
\newblock Differential evolution --- a simple and efficient heuristic for
  global optimization over continuous spaces.
\newblock \emph{Journal of Global Optimization}, 11\penalty0 (4):\penalty0
  341--359, 1997.

\end{thebibliography}

\end{document}